\documentclass[10pt]{article}

\usepackage{a4wide}
\usepackage{amssymb}
\usepackage{amsfonts}
\usepackage{amsmath}
\input xy
\xyoption{arrow} \xyoption{matrix}

\date{}

\newtheorem{proposition}{Proposition}[section]
\newtheorem{theorem}[proposition]{Theorem}
\newtheorem{lemma}[proposition]{Lemma}

\def\der{\partial }

\def\nFM0{{\nu }_{F,M_0}}
\def\nFN0{{\nu }_{F,N_0}}
\def\nGN0{{\nu }_{G,N_0}}

\def\N0{ {\bf N}_0 }

\def\t{\otimes}

\def\v{\varphi}
\def\ra{\rightarrow}

\def\Xpm{X^{\pm }}

\def\s{\sigma}
\def\Z{\mathbb{Z}}

\def\l1{{\lambda}_1}

\def\a{\alpha}
\def\a0{ {\alpha }_0}
\def\a1{ {\alpha }_1}

\def\l{\lambda}

\def\nFGM0{{\nu }_{F,G,M_0}}

\def\nFN0{{\nu}_{F,N_0}}

\def\sm{{\sigma}^m}

\def\sm1{{\sigma}^{-1}}

\def\smtp1{{\sigma}^{-t+1}}

\def\S1{S^{-1}}

\def\Xpm1{X^{\pm 1}_1}

\def\sPM1{{\sigma }^{\pm 1}}
\def\sMP1{{\sigma }^{\mp 1 }}

\def\d{\delta}

\def\di{{\rm d.ind}}

\def\L{\Lambda}

\def\Ytm1{Y^{t-1}}
\def\Yim1{Y^{i-1}}

\def\CK{{\cal K}}

\def\Aut{{\rm Aut}}

\def\dim{{\rm dim }}

\def\ker{ {\rm ker } }

\def\SL2Z{ {\rm SL}_2({\bf Z}) }

\def\Gp1{ G^{1 , 1 } }
\def\P11{ P^{-1 , 1 } }
\def\Pp1{ P^{1 , 1 } }

\def\nCLsr{{}^\nu\kern-2pt {\cal L}^{\sigma , \rho  }}
\def\nP{{}^\nu \kern-2pt P}
\def\nL{{}^\nu\kern-2pt L}
\def\nLL{{}^\nu\kern-2pt \Lambda}
\def\nPsr{{}^\nu\kern-2pt P^{\sigma , \rho  }}
\def\nLsr{{}^\nu\kern-2pt L^{\sigma , \rho  }}
\def\nuCL{{}^\nu\kern-2pt  {\cal L}}
\def\nCLsr{{}^\nu\kern-2pt {\cal L}^{\sigma , \rho  }}
\def\nCL1m{{}^\nu\kern-2pt {\cal L}^{-1 , 1  }}
\def\x1nu{x^\frac{1}{\nu}}
\def\xm1nu{x^{-\frac{1}{\nu}}}

\def\ra{\rightarrow }

\def\CB{{\cal B}}

\def\CC{ {\cal C}}

\def\nAM0{{\nu }_{{\cal A},M_0}}
\def\nAN0{{\nu }_{{\cal A},N_0}}

\def\End{ {\rm End }}

\def\ga{\mathfrak{a}}

\def\SL{{\rm SL}}

\def\di!{\frac{\der^i}{i!}}
\def\dik!{\frac{\der^k_i}{k!}}

\def\N{\mathbb{N}}

\def\0{\overline{0}}
\def\1{\overline{1}}

\def\Ln1{\L_{n,\overline{1}}}

\def\oa{\overline{a}}

\def\a1{a_{\overline{1}}}

\def\S{\Sigma}

\def\vn1{\overrightarrow{n-1}}

\def\soc{{\rm soc}}

\def\mJ{\mathbb{J}}
\def\mI{\mathbb{I}}

\def\K1{{\rm K}_1}

\def\hmI1{\widehat{\mI_1}}
\def\tmI1{\widetilde{\mI_1}}
\def\tmJ1{\widetilde{\mJ_1}}
\def\hB1{\widehat{B_1}}
\def\hCB1{\widehat{\CB_1}}

\def\lsoc{{\rm l.soc}}

\begin{document}

\author{V. V. \  Bavula  
}

\title{The  algebra of polynomial integro-differential operators
 is a holonomic bimodule over the subalgebra of polynomial differential operators}

\maketitle

\begin{abstract}
In contrast to its subalgebra $A_n:=K\langle x_1, \ldots , x_n,
\frac{\der}{\der x_1}, \ldots ,\frac{\der}{\der x_n}\rangle $ of
polynomial  differential operators (i.e. the $n$'th Weyl algebra),
the algebra    $\mI_n:=K\langle x_1, \ldots , x_n,
\frac{\der}{\der x_1}, \ldots ,\frac{\der}{\der x_n}, \int_1,
\ldots , \int_n\rangle $ of polynomial  integro-differential
operators is neither left nor right Noetherian algebra; moreover
it contains infinite direct sums of nonzero left and right ideals.
It is proved that $\mI_n$ is a left (right) coherent algebra iff
$n=1$;   the  algebra $\mI_n$  is a {\em holonomic $A_n$-bimodule}
of length $3^n$ and has multiplicity $3^n$, and all $3^n$ simple
factors of $\mI_n$ are pairwise non-isomorphic $A_n$-bimodules.
 The socle length of the  $A_n$-bimodule
$\mI_n$ is $n+1$, the socle filtration is found, and the $m$'th
term of the socle filtration has length ${n\choose m}2^{n-m}$.
This fact gives a new canonical form for each polynomial
integro-differential operator. It is proved that the algebra
$\mI_n$ is the maximal left (resp. right) order in the largest
left (resp. right) quotient ring of the algebra $\mI_n$.




\end{abstract}


\section{Introduction}

Throughout, ring means an associative ring with $1$; module means
a left module;
 $\N :=\{0, 1, \ldots \}$ is the set of natural numbers; $K$ is a
field of characteristic zero and  $K^*$ is its group of units;
$P_n:= K[x_1, \ldots , x_n]$ is a polynomial algebra over $K$;
$\der_1:=\frac{\der}{\der x_1}, \ldots , \der_n:=\frac{\der}{\der
x_n}$ are the partial derivatives ($K$-linear derivations) of
$P_n$; $\End_K(P_n)$ is the algebra of all $K$-linear maps from
$P_n$ to $P_n$;
the subalgebras   $A_n:= K \langle x_1, \ldots , x_n , \der_1,
\ldots , \der_n\rangle$  and $\mI_n:=K\langle x_1, \ldots , x_n$,
 $\der_1, \ldots ,\der_n,  \int_1,
\ldots , \int_n\rangle $  of the algebra $\End_K(P_n)$ are called
the $n$'th {\em Weyl} algebra and the {\em algebra of  polynomial
integro-differential operators} respectively.

The Weyl algebras $A_n$ are Noetherian algebras and domains. The
algebras $\mI_n$ are neither left nor right Noetherian and not
domains. Moreover, they contain infinite direct sums of nonzero
left and right ideals \cite{Bav-algintdif}. The algebra $A_n$ is
isomorphic to its opposite algebra $A_n^{op}$ via the $K$-algebra
{\em involution}:
$$ A_n\ra A_n, \;\; x_i\mapsto \der_i, \;\; \der_i\mapsto x_i,
\;\; i=1, \ldots , n.$$ Therefore, every $A_n$-bimodule is a left
$A_{2n}$-module and vice versa. {\em Inequality of Bernstein}
\cite{Bernstein-1972} states that {\em each nonzero finitely
generated $A_n$-module has Gelfand-Kirillov dimension which is
greater or equal to $n$}. A finitely generated $A_n$-module is
{\em holonomic} if it has Gelfand-Kirillov dimension $n$. The
holonomic $A_n$-modules share many pleasant properties. In
particular, all holonomic modules have finite length, each nonzero
submodule and factor module of a holonomic  module is holonomic.
The aim of the paper is to prove Theorem \ref{15Sep10}. In
particular, to show that the algebra $\mI_n$ is a holonomic
$A_n$-bimodule of length $3^n$ and has multiplicity $3^n$, i.e. a
holonomic left $A_{2n}$-module of length $3^n$ and has
multiplicity $3^n$. All $3^n$ simple factors of $\mI_n$ are
pairwise non-isomorphic $A_n$-bimodules.  We also found the socle
filtration of the $A_{2n}$-module $\mI_n$. It turns out that the
socle length of the $A_{2n}$-module is $n+1$, and the length, as
an $A_{2n}$-module, of the $m$'th socle factor is ${n\choose m}
2^{n-m}$ (Theorem \ref{15Sep10}.(4)) where $m=0, 1, \ldots , n$. A
new $K$-basis for the algebra $\mI_n$ is found which gives a new
canonical form for each polynomial integro-differential operator,
see (\ref{AFIs1}). By the very definition,
$\mI_n=\bigotimes_{i=1}^n\mI_1(i)\simeq \mI_1^{\t n}$ where
$\mI_1(i):= K\langle x_i, \der_i, \int_i\rangle$ and
$A_n=\bigotimes_{i=1}^n A_1(i) = A_1^{\t n }$ where $A_1(i) :=
K\langle x_i ,\der_i \rangle$. So, the properties of the algebras
$\mI_n$ and $A_n$ are `determined' by the properties of the
algebras $\mI_1$ and $A_1$.

At the beginning of Section \ref{PTH15} we collect necessary facts
on the algebras $\mI_n$. Then we prove Theorem \ref{15Sep10} in
the case when $n=1$ and prove some necessary results that are used
in the proof of Theorem \ref{15Sep10} (in the  general case) which
is  given at the end of the section.

In Section \ref{AInCOH}, it is proved that the algebra $\mI_n$ is
left (right) coherent iff $n=1$ (Theorem \ref{30Oct10}).

In Section \ref{AIMAXO}, it is proved that the algebra $\mI_n$ is
the maximal left (resp. right) order in its largest left (resp.
right) quotient ring (Theorem \ref{14Feb11}).


\section{Proof of Theorem \ref{15Sep10}}\label{PTH15}

At the beginning of this section, we collect necessary (mostly
elementary) facts on the algebra $\mI_1$ from \cite{Bav-algintdif}
that are used later in the paper.

The algebra $\mI_1$  is generated by the elements $\der $, $H:=
\der x$ and $\int$ (since $x=\int H$) that satisfy the defining
relations (Proposition 2.2, \cite{Bav-algintdif}):
\begin{equation}\label{I1rel}
\der \int = 1, \;\; [H, \int ] = \int, \;\; [H, \der ] =-\der ,
\;\; H(1-\int\der ) =(1-\int\der ) H = 1-\int\der ,
\end{equation}
 where
$[a,b]:=ab-ba$ is the {\em commutator} of elements $a$ and $b$.
 The elements of the algebra $\mI_1$,  
\begin{equation}\label{eijdef}
e_{ij}:=\int^i\der^j-\int^{i+1}\der^{j+1}, \;\; i,j\in \N ,
\end{equation}
satisfy the relations $e_{ij}e_{kl}=\d_{jk}e_{il}$ where $\d_{jk}$
is the Kronecker delta function and $\N := \{ 0, 1, \ldots \}$ is
the set of natural numbers. Notice that
$e_{ij}=\int^ie_{00}\der^j$. The matrices of the linear maps
$e_{ij}\in \End_K(K[x])$ with respect to the basis $\{ x^{[s]}:=
\frac{x^s}{s!}\}_{s\in \N}$ of the polynomial algebra $K[x]$  are
the elementary matrices, i.e.
$$ e_{ij}*x^{[s]}=\begin{cases}
x^{[i]}& \text{if }j=s,\\
0& \text{if }j\neq s.\\
\end{cases}$$
Let $E_{ij}\in \End_K(K[x])$ be the usual matrix units, i.e.
$E_{ij}*x^s= \d_{js}x^i$ for all $i,j,s\in \N$. Then
\begin{equation}\label{eijEij}
e_{ij}=\frac{j!}{i!}E_{ij},
\end{equation}
 $Ke_{ij}=KE_{ij}$, and
$F:=\bigoplus_{i,j\geq 0}Ke_{ij}= \bigoplus_{i,j\geq
0}KE_{ij}\simeq M_\infty (K)$, the algebra (without 1) of infinite
dimensional matrices. $F$ is the only proper ideal (i.e. $\neq 0,
\mI_1$) of the algebra $\mI_1$ \cite{Bav-algintdif}.

$\noindent $

{\bf $\Z$-grading on the algebra $\mI_1$ and the canonical form of
an integro-differential operator \cite{Bav-algintdif},
\cite{Bav-intdifline}}. The algebra $\mI_1=\bigoplus_{i\in \Z}
\mI_{1, i}$ is a $\Z$-graded algebra ($\mI_{1, i} \mI_{1,
j}\subseteq \mI_{1, i+j}$ for all $i,j\in \Z$) where
$$ \mI_{1, i} =\begin{cases}
D_1\int^i=\int^iD_1& \text{if } i>0,\\
D_1& \text{if }i=0,\\
\der^{|i|}D_1=D_1\der^{|i|}& \text{if }i<0,\\
\end{cases}
 $$
 the algebra $D_1:= K[H]\bigoplus \bigoplus_{i\in \N} Ke_{ii}$ is
a {\em commutative non-Noetherian} subalgebra of $\mI_1$, $
He_{ii} = e_{ii}H= (i+1)e_{ii}$  for $i\in \N $ (and so
$\bigoplus_{i\in \N} Ke_{ii}$ is the direct sum of non-zero ideals
$Ke_{ii}$ of the algebra $D_1$); $(\int^iD_1)_{D_1}\simeq D_1$,
$\int^id\mapsto d$; ${}_{D_1}(D_1\der^i) \simeq D_1$,
$d\der^i\mapsto d$,   for all $i\geq 0$ since $\der^i\int^i=1$.
 Notice that the maps $\cdot\int^i : D_1\ra D_1\int^i$, $d\mapsto
d\int^i$,  and $\der^i \cdot : D_1\ra \der^iD_1$, $d\mapsto
\der^id$, have the same kernel $\bigoplus_{j=0}^{i-1}Ke_{jj}$.

Each element $a$ of the algebra $\mI_1$ is the unique finite sum
\begin{equation}\label{acan}
a=\sum_{i>0} a_{-i}\der^i+a_0+\sum_{i>0}\int^ia_i +\sum_{i,j\in
\N} \l_{ij} e_{ij}
\end{equation}
where $a_k\in K[H]$ and $\l_{ij}\in K$. This is the {\em canonical
form} of the polynomial integro-differential operator
\cite{Bav-algintdif}.

$\noindent $

{\it Definition}. Let $a\in \mI_1$ be as in (\ref{acan}) and let
$a_F:=\sum \l_{ij}e_{ij}$. Suppose that $a_F\neq 0$ then
\begin{equation}\label{degFa}
\deg_F(a) :=\min \{ n\in \N \, | \, a_F\in \bigoplus_{i,j=0}^n
Ke_{ij}\}
\end{equation}
is called the $F$-{\em degree} of the element $a$;
$\deg_F(0):=-1$.

$\noindent $

Let $v_i:=\begin{cases}
\int^i& \text{if }i>0,\\
1& \text{if }i=0,\\
\der^{|i|}& \text{if }i<0.\\
\end{cases}$
Then $\mI_{1,i}=D_1v_i= v_iD_1$ and an element $a\in \mI_1$ is the
unique  finite  sum 
\begin{equation}\label{acan1}
a=\sum_{i\in \Z} b_iv_i +\sum_{i,j\in \N} \l_{ij} e_{ij}
\end{equation}
where $b_i\in K[H]$ and $\l_{ij}\in K$. So, the set $\{ H^j\der^i,
H^j, \int^iH^j, e_{st}\, | \, i\geq 1; j,s,t\geq 0\}$ is a
$K$-basis for the algebra $\mI_1$. The multiplication in the
algebra $\mI_1$ is given by the rule:
$$ \int H = (H-1) \int , \;\; H\der = \der (H-1), \;\; \int e_{ij}
= e_{i+1, j}, \;\; e_{ij}\int= e_{i,j-1}, \;\; \der e_{ij}=
e_{i-1, j}\;\; e_{ij} \der = \der e_{i, j+1}.$$
$$ He_{ii} = e_{ii}H= (i+1)e_{ii}, \;\; i\in \N, $$
where $e_{-1, j}:=0$ and $e_{i,-1}:=0$.

 $\noindent $

The factor algebra $B_1:= \mI_1/F$ is the simple Laurent skew
polynomial algebra $K[H][\der, \der^{-1}; \tau ]$ where the
automorphism $\tau \in \Aut_{K-{\rm alg}}(K[H])$ is defined by the
rule $\tau (H) = H+1$, \cite{Bav-algintdif}. Let
\begin{equation}\label{piI1B1}
\pi : \mI_1\ra B_1, \;\; a\mapsto \oa : a+F,
\end{equation}
be the canonical epimorphism.

$\noindent $

The Weyl algebra $A_2$ is equipped with the, so-called, {\em
filtration of Bernstein}, $A_2=\bigcup_{i\geq 0}A_{2, \leq i}$
where $A_{2, \leq i}:= \bigoplus \{
Kx_1^{\alpha_1}x_2^{\alpha_2}\der_1^{\beta_1}\der_2^{\beta_2}\, |
\, \alpha_1+\alpha_2+\beta_1+\beta_2\leq i\}$. The polynomial
algebra $P_2:= K[x_1, x_2]\simeq A_2/(A_2\der_1+A_2\der_2)$ is a
simple left $A_2$-module with $\End_{A_2} (P_2)=
\ker_{P_2}(\der_1) \cap \ker_{P_2}(\der_2) = K$. The standard
filtration $\{ A_{2, \leq i}\cdot 1\}_{i\in \N}$ of the
$A_2$-module $P_2$ coincides with the filtration $\{ P_{2, \leq
i}:=\sum_{\alpha_1, \alpha_2\geq 0} \{
Kx_1^{\alpha_1}x_2^{\alpha_2}\, | \, \alpha_1+\alpha_2\leq i
\}_{i\in \N}$ on the polynomial algebra $P_2$ by the total degree,
i.e. $P_{2, \leq i}=A_{2, \leq i}\cdot 1$ for all $i\geq 0$, and
so $\dim_K(A_{2, \leq i} ) = \frac{(i+1)(i+2)}{2}$. Therefore,
$P_2$ is a holonomic $A_2$-module with multiplicity $e(P_2)=1$ and
$\End_{A_2}(P_2)\simeq K$. The Weyl algebra $A_1$ admits the
$K$-isomorphism: 
\begin{equation}\label{A1xd}
\xi : A_1, \ra A_1, \;\; x\mapsto \der , \;\; \der\mapsto -x.
\end{equation}
Then $1\t \xi$ is an automorphism of the Weyl algebra $A_2$. The
twisted by the automorphism $1\t \xi$ $A_2$-module $P_2$,
\begin{equation}\label{P21xi}
P_2^{1\t \xi} \simeq K[x_1, \der_2] \simeq A_2/(A_2\der_1+A_2x_2)
\end{equation}
is a simple holonomic $A_2$-module with multiplicity 1 and
$\End_{A_2}(P_2^{1\t \xi})\simeq K$.

The Weyl algebra $A_1$ is isomorphic to its opposite algebra
$A_1^{op}$ via 
\begin{equation}\label{A1op}
A_1\ra A_1^{op}, \;\; x\mapsto \der , \;\; \der \mapsto x.
\end{equation}
In particular, each $A_1$-bimodule ${}_{A_1}M_{A_1}$ is a left
$A_2$-module: ${}_{A_1}M_{A_1}={}_{A_1\t A_1^{op}}M\simeq
{}_{A_1\t A_1}M= {}_{A_2}M$.

\begin{lemma}\label{a12Sep10}
\begin{enumerate}
\item ${}_{A_1}F_{A_1}= A_1e_{00}A_1\simeq {}_{A_1}(A_1/A_1\der \t
A_1/xA_1)_{A_1}$. \item ${}_{A_2}F\simeq A_2/(A_2\der_1+A_2\der_2
)\simeq K[x_1, x_2]$ is a simple holonomic $A_2$-module with
multiplicity 1 with respect to the filtration of Bernstein of the
algebra  $A_2$ and $End_{A_2}(F)\simeq K$.
\end{enumerate}
\end{lemma}

{\it Proof}.   ${}_{A_1}(A_1/A_1\der \t A_1 / xA_1)_{A_1} \simeq
{}_{A_1\t A_1}(A_1/A_1\der \t A_1 / A_1\der )\simeq
A_2/(A_2\der_1+A_2\der_2 )\simeq K[x_1, x_2]$ is a simple
holonomic $A_2$-module with multiplicity 1  with respect to the
filtration of Bernstein of the algebra  $A_2$ and
$\End_{A_2}(F)\simeq K$. Clearly, ${}_{A_1}F_{A_1}= A_1e_{00}A_1$
 and the $A_1$-bimodule homomorphism
 $$A_1/A_1\der \t A_1 / xA_1\ra A_1e_{00}A_1, \;\; (1+A_1\der_1)\t
 (1+xA_1)\mapsto e_{00},$$
 is an epimorphism. Therefore, it is an isomorphism by the
 simplicity of the first $A_1$-bimodule.
$\Box $


\begin{proposition}\label{b12Sep10}
\begin{enumerate}
\item ${}_{A_1}(\mI_1 / (A_1+F))_{A_1}\simeq A_1/A_1\der \t A_1 /
\der A_1$. \item ${}_{A_2}(\mI_1/(A_1+F))\simeq A_1/A_1\der \t A_1
/ A_1x\simeq A_2/(A_2\der_1+A_2x_2)\simeq K[x_1, \der_2]$ is a
simple holonomic $A_2$-module with multiplicity 1 with respect to
the filtration of Bernstein and $\End_{A_2}(K[x_1, \der_2])\simeq
 K$. \item ${}_{A_1}(\mI_1/(A_1+F))\simeq (A_1/A_1\der )^{(\N )}
\simeq K[x]^{(\N )}$ is a semi-simple left $A_1$-module and
$(\mI_1/(A_1+F))_{A_1} \simeq (A_1/ \der A_1)^{(\N )}\simeq
K[x]^{(\N )}$ is a semi-simple right $A_1$-module.
\end{enumerate}
\end{proposition}

{\it Proof}. 1 {\em and } 2.  Notice that ${}_{A_2}(A_1/A_1\der \t
A_1 / A_1x)\simeq {}_{A_2}(A_2/(A_2\der_1+A_2x_2))\simeq K[x_1,
\der_2]$ is a simple holonomic $A_2$-module with multiplicity 1
with respect to the filtration of Bernstein  and
$\End_{A_2}(K[x_1, \der_2])\simeq K$.  The natural filtration of
the polynomial algebra $Q':=K[x_1, \der_2]=\bigcup_{i\geq 0}
Q'_{\leq i}$ by the total degree of the variables, i.e. $Q_{\leq
i}':= \bigoplus_{s+t\leq i} Kx_1^s\der_2^t$,  is a standard
filtration for the $A_2$-module $Q'=A_2\cdot 1$ since $Q'_{\leq
i}= A_{2, \leq i}\cdot 1$ for all $i\geq 0$. In particular,
$\dim_K(Q'_{\leq i} )= \frac{(i+1)(i+2)}{2}$ for all $i\geq 0$. By
(\ref{acan}), the $A_1$-bimodule $Q:= \mI_1 / (A_1+F)$ is the
direct sum 
\begin{equation}\label{QIQi}
Q=\bigoplus_{i\geq 1}Q_i
\end{equation}
of its vector subspaces 
\begin{equation}\label{1QIQi}
({Q_i})_{K[H]}\simeq \int^i K[H]/x^iK[H]\simeq \int^i K[H]/\int^i
(H(H+1)\cdots (H+i-1))\simeq K[H] / (H(H+1)\cdots (H+i-1))
\end{equation}
(since $x^i=(\int H)^i = \int^i H(H+1)\cdots (H+i-1)$ and $ \der^i
\int^i =1$) such that $xQ_i\subseteq Q_{i+1}$, $Q_ix \subseteq
Q_{i+1}$, $\der Q_i \subseteq Q_{i-1}$ and $Q_i\der \subseteq
Q_{i-1}$ for all all $i\geq 1$ where $Q_{0}:=0$.  Then
$A_2$-module $Q$ has the finite dimensional ascending filtration
$Q= \bigcup_{i\geq 0} Q_{\leq i} $ where $Q_{\leq i}:=
\bigoplus_{1\leq j\leq i+1}Q_j$ and
$$\dim_K(Q_{\leq i}) =
\sum_{j=0}^i (j+1) = \frac{(i+1)(i+2)}{2}\;\; {\rm  for\; all}\;\;
i\geq 0.$$ Since $\der Q_1 = Q_1\der =0$, the simple filtered
$A_2$-module (treated as $A_1$-bimodule) ${}_{A_1}Q'_{A_1}=
A_1/A_1\der \t A_1 /\der A_1$ can be seen as a filtered
$A_2$-submodule of $Q$ via $(1+A_1\der ) \t (1+ \der A_1)\mapsto
\int +A_1+F$. In particular, for all $i\geq 0$, we have the
inclusions $Q_i'\subseteq Q_i$ which are, in fact, equalities
since $\dim_K(Q_i')= \dim_K(Q_i)$. Then,
$$ {}_{A_1}(\mI_1/ (A_1+F))_{A_1} \simeq {}_{A_1}(A_1/A_1\der \t
A_1/\der A_1)_{A_1} \simeq {}_{A_2}(A_1/A_1\der \t A_1 / A_1
x)\simeq K[x_1, \der_2].$$ It is obvious that the $A_2$-module
$K[x_1, \der_2]$ is a simple $A_2$-module with multiplicity 1 and
$\End_{A_2} (K[x_1, \der_2])$ $ \simeq K$.

3. Statement 3 follows from statement 1.  $\Box $

$\noindent $

A linear map $\v$ acting in a vector space $V$ is called a {\em
locally nilpotent map} if $V=\bigcup_{i\geq 1} \ker (\v^i)$, i.e.
for each element $v\in V$ there exists a natural number $i$ such
that $\v^i (v)=0$.

$\noindent $

It follows from Proposition \ref{b12Sep10} and (\ref{1QIQi}) that
\begin{equation}\label{kerdAF}
\ker_{\mI_1/ (A_1+F)}(\der \cdot) \bigcap \ker_{\mI_1/
(A_1+F)}(\cdot \der )= K (\int +A_1+F),
\end{equation}
and that the maps $\der\cdot : \mI_1/ (A_1+F)\ra \mI_1/ (A_1+F)$,
$u\mapsto \der u$,  and $\cdot \der: \mI_1/ (A_1+F)\ra \mI_1/
(A_1+F)$, $u\mapsto u\der $, are locally nilpotent since
\begin{equation}\label{dxx}
\der*x_1^ix_2^j= ix_1^{i-1}x_2^j, \;\;  x_1^ix_2^j*\der =-
jx_1^ix_2^{j-1}.
\end{equation}
   Recall that the
{\em socle} $\soc_A(M)$ of a module $M$ over a ring $A$ is the sum
of all the simple submodules of $M$, if they exist, and zero,
otherwise.

\begin{theorem}\label{13Sep10}
\begin{enumerate}
\item The $A_1$-bimodule $\mI_1$ is a holonomic $A_2$-module of
length 3 with simple non-isomorphic factors $F\simeq
{}_{A_2}K[x_1, x_2]$, ${}_{A_1}{A_1}_{A_1}$ and ${}_{A_2}K[x_1,
\der_2]$. Each factor is a simple holonomic $A_2$-module with
multiplicity 1 and its $A_2$-module endomorphism algebra  is $K$.
\item $\soc_{A_2}(\mI_1)=A_1\bigoplus F$.  \item The short exact
sequence of $A_2$-modules 
\begin{equation}\label{AFIs}
0\ra A_1\bigoplus F\ra \mI_1 \ra \mI_1/ (A_1+F)\ra 0
\end{equation}
is non-split.
\end{enumerate}
\end{theorem}

{\it Proof}. 1. Statement 1 follows from Lemma \ref{a12Sep10},
Proposition \ref{b12Sep10} and (\ref{AFIs}).

3. Suppose that the short exact sequence of $A_1$-bimodules
splits, we seek a contradiction. Then, by Proposition
\ref{b12Sep10}.(1) and (\ref{kerdAF}), there is a nonzero element,
say $u=\int +a+f\in \mI_1$ with $a\in A_1$ and $f\in F$ such that
$\der u =0$ and $u\der =0$. The first equation implies $1+ \der
a=-\der f\in A_1 \cap F=0$, and so $\der a = -1$ in $A_1$, a
contradiction.

2. Statement 2 follows from statement 3. $\Box $

$\noindent $

{\bf New basis for the algebra $\mI_n$}. It follows from
(\ref{QIQi}), (\ref{1QIQi}) and (\ref{AFIs}) that
\begin{equation}\label{AFIs1}
\mI_1=\bigoplus_{i,j\geq 0}Kx^i\der^j\oplus \bigoplus_{k,l\geq 0}
Ke_{kl}\oplus \bigoplus\{ K\int^sH^t\, | \, s\geq 1, t=0, 1,
\ldots , s-1\} .
\end{equation}
This gives a new $K$-basis for the algebra $\mI_1$: $\{ x^i\der^j,
e_{kl}, \int^sH^t\, | \, i,j,k,l\geq 0; s\geq 1; t=0, 1, \ldots ,
s-1\}$. By taking $n$'th tensor product of this basis we obtain a
new $K$-basis for the algebra $\mI_n = \mI_1^{\t n}$.

\begin{lemma}\label{a14Sep10}
\begin{enumerate}
\item The $A_1$-bimodule $\mI_1/ A_1$ is a holonomic $A_2$-module
of length 2 with simple non-isomorphic factors $F\simeq {}_{A_2}
K[x_1, x_2]$ and ${}_{A_2} K[x_1, \der_2]$. Each factor is a
simple holonomic $A_2$-module with multiplicity 1 and its
$A_2$-module endomorphism algebra is $K$. \item
$\soc_{A_2}(\mI_1/A_1)=F$. \item The short exact sequence of
$A_2$-modules 
\begin{equation}\label{FIA}
0\ra F\ra \mI_1/A_1\ra \mI_1/(A_1+F)\ra 0
\end{equation}
is non-split. \item The short exact sequence of left $A_1$-modules
(\ref{FIA}) splits and so ${}_{A_1}(\mI_1/A_1)\simeq K[x]^{(\N )}$
is a semi-simple left $A_1$-module. \item The short exact sequence
of right $A_1$-modules (\ref{FIA}) does not split, and so $(\mI_1
/ A_1)_{A_1}$ is not a semi-simple right $A_1$-module.
\end{enumerate}
\end{lemma}

{\it Proof}. 1. Statement 1 follows from Theorem
\ref{13Sep10}.(1).

3. Suppose that the short exact sequence of $A_1$-bimodules
(\ref{FIA}) splits, we seek a contradiction. Then, by Proposition
\ref{b12Sep10}.(1) and (\ref{kerdAF}), there is a nonzero element,
say $u=\int +f+A_1\in \mI_1 / A_1$ with $f\in F$ such that $0=\der
u = 1+\der f$ and $ 0 = u\der = 1-e_{00}+f\der$ in $\mI_1/ A_1$.
The first equality gives $\der f =0$ in $\mI_1/ A_1$, and so $f =
\sum_{i\geq 0} \l_i e_{0i}$ for some $\l_i\in K$. Then the second
equality gives $e_{00}=f\der  = \sum_{i\geq 0} \l_i e_{0i}\der =
\sum_{i\geq 0} \l_i e_{0,i+1}$, a contradiction.

2. Statement 2 follows from statement 3.

4.  Let $L$ be the last sum in the decomposition (\ref{AFIs1}),
i.e. 
\begin{equation}\label{AFIs2}
\mI_1= A_1\bigoplus F\bigoplus L.
\end{equation}
Then $A_1\bigoplus L$ is a left $A_1$-submodule of ${}_{A_1}\mI_1$
since $\der \int =1$, $x=\int H$ and $\int H = (H-1) \int$. Notice
that $A_1\bigoplus L$ is not a right $A_1$-submodule of $\mI_1$
since $\int\der = 1-e_{00}\not\in A_1\bigoplus L$.
 By (\ref{AFIs2}), ${}_{A_1}(\mI_1 / A_1) \simeq F\bigoplus
 (A_1+L)/A_1$ is a direct sum of left $A_1$-submodules such that
 ${}_{A_1}F\simeq K[x]^{(\N )}$ (Lemma \ref{a12Sep10}.(1)) and
 ${}_{A_1}((A_1+L)/A_1) \simeq \mI_1 / (A_1+F)\simeq K[x]^{(\N )}$
  (Proposition \ref{b12Sep10}.(3)). Therefore,
  ${}_{A_1}(\mI_1/A_1)$ is a semi-simple module.  Therefore,  the short exact
  sequence of left $A_1$-modules (\ref{FIA}) splits and
  ${}_{A_1}(\mI_1/ A_1)\simeq K[x]^{(\N )}\bigoplus K[x]^{(\N
  )}\simeq K[x]^{(\N )}$.

  5. By Proposition \ref{b12Sep10}.(1), $(\mI_1/
  (A_1+F))_{A_1}\simeq (A_1 / \der A_1)^{(\N )}$. Suppose that the
  short exact sequence of right $A_1$-modules (\ref{FIA}) splits,
  we seek a contradiction. In the factor module $\mI_1 / (A_1+F)$,
  $(\int +A_1+F) \der =0$ since $\int\der = 1-e_{00}\in A_1+F$.
  Then the splitness   implies that $(\int +f+A_1)\der =0$ in $\mI_1
  / A_1$ for some element $f\in F$, equivalently, $-e_{00}+f\der \in
  A_1\cap F=0$ in $\mI_1$, i.e. $f\der = e_{00}$, this is obviously
  impossible (since $e_{i,j}\der = e_{i,j+1}$), a contradiction. $\Box $

$\noindent $

Let $M$ be a module over a ring $R$. The socle $\soc_R(M)$, if
nonzero, is the largest semi-simple submodule of $M$. The socle of
$M$, if nonzero, is the only  essential  semi-simple  submodule.
The {\em socle chain} of the module $M$ is the ascending chain of
its submodules:
$$ \soc_R^0(M):=\soc_R(M) \subseteq \soc_R^1(M) \subseteq\cdots
\subseteq \soc_R^i(M) \subseteq\cdots $$ where $\soc_R^i(M):=
\v_{i-1}^{-1} (\soc_R (M/\soc_R^{i-1}(M)))$ where $\v_{i-1}:M\ra
M/\soc^{i-1}(M)$, $m\mapsto m+\soc_R^{i-1}(M)$. Let $\soc_R^\infty
(M) := \bigcup_{i\geq 0}\soc_R^i (M)$. If $M=\soc_R^\infty (M)$
then
$$ \lsoc_R(M) =1+\min \{ i\geq 0\, | \, M= \soc_R^i(M)\}$$
is called the {\em socle length} of the $R$-module $M$. So, a
nonzero module is semi-simple iff its socle length is 1.

\begin{theorem}\label{15Sep10}
\begin{enumerate}
\item The $A_n$-bimodule $\mI_n$ is a holonomic $A_{2n}$-module of
length $3^n$ with pairwise non-isomorphic simple factors and each
of them  is the tensor product $\bigotimes_{i=1}^n M_i$  of simple
$A_2(i)$-modules $M_i$ as in  Theorem \ref{13Sep10} for $i=1,
\ldots , n$. Each simple factor $\bigotimes_{i=1}^n M_i$ is a
simple holonomic $A_{2n}$-module and has multiplicity 1 (with
respect to the filtration of Bernstein on the algebra $A_{2n}$)
and its $A_{2n}$-module endomorphism algebra is $K$.
 \item $\soc_{A_{2n}}(\mI_n) = \bigotimes_{i=1}^n \soc_{A_2(i)}(\mI_1(i))=
 \bigotimes_{i=1}^n (A_1(i) \bigoplus F(i))$.
 \item The socle length of the $A_{2n}$-module $\mI_n$ is $n+1$.
 For each number $m=0,1, \ldots , n$, $$\soc_{A_{2n}}^m(\mI_n) =
 \sum_{i_1+\cdots +i_n=m} \bigotimes_{s=1}^n \soc_{A_2(i)}^{i_s}
 (\mI_1(i))$$  where all $i_s \in \{ 0,1\}$ and $\soc^j_{A_2(i)}
 = \begin{cases}
A_1(i) \bigoplus F(i)& \text{if }j=0,\\
\mI_1(i) & \text{if }j=1.\\
\end{cases}$
\item For each number $m=0,1, \ldots , n$,
$$ \soc_{A_{2n}}^m (\mI_n)/ \soc_{A_{2n}}^{m-1} (\mI_n)=
\bigoplus_{i_1+\cdots +i_n=m} \bigotimes_{s=1}^n
\soc_{A_2(i)}^{i_s}
 (\mI_1(i))/\soc_{A_2(i)}^{i_s-1}
 (\mI_1(i))$$
and its length (as an $A_{2n}$-module) is ${n\choose m} 2^{n-m}$
where all $i_s \in \{ 0,1\}$ and $\soc^{-1}:=0$. \item The left
$A_{2n}$-module $\mI_n$ has multiplicity $3^n$ with respect to the
filtration of Bernstein of the Weyl algebra $A_{2n}$.
\end{enumerate}
\end{theorem}

{\em Remark}. The sum of lengths of all the factors in statement 4
is $3^n$ as $$3^n = (1+2)^n = \sum_{m=0}^n {n\choose m} 2^{n-m}.$$

{\it Proof}. 1. By Theorem \ref{13Sep10}.(1), each of the tensor
multiples $\mI_1(i)$ in $\mI_n =\bigotimes_{i=1}^n \mI_1(i)$ has
the $A_2(i)$-module (i.e. the $A_1(i)$-bimodule) filtration of
length 3 with factors $M_i$ as in Theorem \ref{13Sep10}.(1). By
considering the tensor product of these filtrations, the
$A_{2n}$-module $\mI_n = \bigotimes_{i=1}^n \mI_1 (i)$ (i.e. the
$A_n$-bimodule) has a filtration (of length $3^n$) with factors
$\bigotimes_{i=1}^n M_i$. It is obvious that each $A_{2n}$-module
$\bigotimes_{i=1}^n M_i$ is isomorphic to a twisted
$A_{2n}$-module ${}^\s P_{2n}$ by an automorphism $\s $ of the
Weyl algebra $A_{2n}$ that preserves the filtration of Bernstein
on the algebra $A_{2n}$ where $$P_{2n}= K[x_1, \ldots ,
x_{2n}]\simeq A_{2n}/ \sum_{i=1}^{2n} A_{2n} \der_i.$$ This
statement is obvious for $n=1$, then the general case follows at
once. Since the $A_{2n}$-module $P_{2n}$ is simple, holonomic with
multiplicity 1 (since $e({}^\s P_{2n}) = e(P_{2n}) = 1$) and
$\End_{A_{2n}}(P_{2n})\simeq K$, then so are all the
$A_{2n}$-modules $\bigotimes_{i=1}^n M_i$. This finishes the proof
of statement 1.

2. Statement 2 follows from statement 3.

3. To prove statement 3 we use induction on $n$. The initial step
when $n=1$ is true due to Theorem \ref{13Sep10}.(1). Suppose that
$n>1$ and the statement holds for all $n'<n$.
 Let $\{ s^0= A_1\bigoplus F, s^1=\mI_1\}$ be the socle filtration
 for ${}_{A_1}{\mI_1}_{A_1}$ and let $\{ s^0, s^1, \ldots ,
 s^{n-1}\}$ be the socle filtration for
 ${}_{A_{n-1}}{\mI_{n-1}}_{A_{n-1}}$. We are going to prove that $$\{
 s'^0:=s^0\t s^0, s'^1:=s^0\t s^1+s^1\t s^0, \ldots , s'^{n-1}:= s^0\t
 s^{n-1} +s^1\t s^{n-2}, s'^n := s^1\t s^{n-1}\}$$ is the socle
 filtration for ${}_{A_n}{\mI_n}_{A_n}$. Notice that $A_n = A_1\t
 A_{n-1}$, $\mI_n = \mI_1\t \mI_{n-1}$ and $\{ s^0\t
 s^i\}_{i=0}^{n-1}$ is the socle filtration for
 ${}_{A_{n-1}}(s^0\t \mI_{n-1})_{A_{n-1}}=s^0\t ({}_{A_{n-1}}{\mI_{n-1}}_{A_{n-1}})$ since the
 $\mI_{n-1}$-bimodules $s^0\t s^i/s^{i-1}$ are semi-simple.
 Since, for each number $m=0,1, \ldots , n$, the $A_n$-subbimodule
 $$ \overline{s}'^m:= s'^m/s'^{m-1} = s^0\t (s^m / s^{m-1})
 \bigoplus (s^1/s^0)\t (s^{m-1}/ s^{m-2})\;\;\;
 ({\rm where} \;\; \overline{s}'^0=s^0\t s^0)$$ of $\mI_n'/ s'^{m-1}$ is
 semi-simple, in order to finish the proof of statement 3 it
 suffices to show that $\overline{s}'^m$ is  an  essential
 $A_n$-subbimodule of $\mI_n/s'^{m-1}$. Let $a$ be a nonzero
  element of the $A_n$-bimodule $\mI_n/ s'^{m-1}$.  We have to show that $A_n
 a A_n \cap \overline{s}'^m\neq 0$. If $a\in s^0 \t
 \mI_{n-1}+s'^{m-1}$ then
 $$0\neq \mI_{n-1} a \mI_{n-1} \cap s^0 \t
 (s^m / s^{m-1}) \subseteq \overline{s}'^m$$ (since
 $\{ s^0\t s^i\}^{n-1}_{i=0}$ is
 the socle filtration for ${}_{A_{n-1}} (s^0\t
 \mI_{n-1})_{A_{n-1} }$). If $a\not\in s^0 \t \mI_{n-1}+s'^{m-1}$
 then using the explicit basis $\{ x_1^ix_2^j\}_{i,j\geq 0}$ for
 the $A_1$-bimodule $s^1/s^0$ (Proposition \ref{b12Sep10}.(1)) and the
 action of the element $\der $ on it (see (\ref{dxx})), we can
 find natural numbers, say $k$ and $l$,  such that, by
 (\ref{kerdAF}), the element
 $$a':= \der^k a \der^l = \int\t
 u_1+v_2\t u_2+\cdots +v_s\t u_s,$$ such that  $0\neq u_1\in \mI_{n-1}/
 s^{m-1}$ (in particular, $a'$ is a nonzero element of $\mI_n/
 s'^{m-1}$); $u_2, \ldots , u_s$ are linearly independent elements
 of $\mI_{n-1}$; $v_2, \ldots , v_s$ are linearly independent
 elements of $s^0$. If the elements $u_1, u_2, \ldots , u_s$ are linearly independent  then
 $$ a'':= \der a'= 1\t u_1+(\der v_2)\t u_2+\cdots + (\der v_s) \t
 u_s$$
 is a nonzero element  of $s^0\t \mI_{n-1}$, and so, by the
 previous case $\mI_n a \mI_n \cap \overline{s}'^m\neq 0$.

 If the elements $u_1, u_2, \ldots , u_s$ are linearly dependent
 then $u_1=\sum_{i=2}^s\l_iu_i$ for some elements $\l_i\in K$ not all
 of which are zero ones, say $\l_2\neq 0$. The element $a'$ can be
 written as $a'=(\l_2\int +v_2) \t u_2+\cdots +(\l_s\int +v_s) \t
 u_s$.
 $$ a'':= \der a' =(\l_2 +\der v_2) \t u_2+\cdots +(\l_s +\der v_s) \t
 u_s.$$
  We claim that $a''\neq 0$. Suppose that $a''=0$, we seek a
  contradiction. Then $\l_2+\der v_2=0, \ldots , \l_s+\der v_s=0$
  in $A_1\bigoplus F$ (since the elements $u_2, \ldots , u_s$
  are linearly independent).  The first equality yields $0\neq \l_2=\der
  b$ in the Weyl algebra $A_1$ for some element $b\in A_1$. This
  is clearly impossible. Therefore, $a''$ is a nonzero element of
  $s^0\t \mI_{n-1}$, and so, by the previous case, $\mI_n a\mI_n
  \cap \overline{s}'^m\neq 0$.

4. The equality follows from statement 3. To prove the claim about
the length note that ${n\choose m}$ is the number of vectors
$(i_1, \ldots , i_n)\in \{ 0,1\}^n$ with $i_1+\cdots + i_n = m$;
and for each choice of $(i_1, \ldots , i_n)$ the length of the
$A_{2n}$-module $\bigotimes_{s=1}^n \soc_{A_2(i)}^{i_s}
 (\mI_1(i))/\soc_{A_2(i)}^{i_s-1}
 (\mI_1(i))$ is $2^{n-m}$. Therefore, the length of the
 $A_{2n}$-module $\soc^m_{A_{2n}}(\mI_n)
 /\soc^{m-1}_{A_{2n}}(\mI_n)$ is ${n\choose m}2^{n-m}$.

5. Statement 5 follows from statement 1 and the additivity of the
multiplicity  on the holonomic modules.

 $\Box $



\section{The algebra $\mI_n$ is coherent iff $n=1$}\label{AInCOH}

The aim of this section is to prove Theorem \ref{30Oct10}.


A module $M$ over a ring $R$ is {\em finitely presented} if there
is an exact sequence of modules $R^m\ra R^n\ra M\ra 0$. A finitely
generated module is a {\em coherent} module if every finitely
generated submodule is finitely presented. A ring $R$  is a {\em
left} (resp. {\em right}) {\em coherent ring} if the module
${}_RR$ (resp. $R_R$) is coherent. {\em A ring $R$ is a left
coherent ring iff, for each element $r\in R$, $\ker_R(\cdot r)$ is
a finitely generated left $R$-module and the intersection of two
finitely generated left ideals is finitely generated}, Proposition
13.3, \cite{Stenstrom-RingQuot}. Each left Noetherian ring is left
coherent but not vice versa.

\begin{theorem}\label{30Oct10}
The algebra $\mI_n$ is a left coherent algebra iff the algebra
$\mI_n$ is a right coherent algebra iff $n=1$.
\end{theorem}

{\it Proof}. The first `iff' is obvious since the algebra $\mI_n$
is self-dual \cite{Bav-algintdif}, i.e. is isomorphic to its
opposite algebra $\mI_n^{op}$. If $n=1$ the algebra is a left
coherent algebra \cite{Bav-intdifline}. If $n\geq 2$ then the
algebra $\mI_2$ is not a left coherent algebra since, by Lemma
\ref{a30Oct10}, $\ker_{\mI_n} (\cdot (H_1-H_2)) = \ker_{\mI_2}
(\cdot (H_1-H_2))\t \mI_{n-2}\simeq {}_{\mI_n}(P_2\t
\mI_{n-2})^{(\N )}$ is an infinite direct sum of nonzero
$\mI_n$-modules, hence it is not finitely generated. Therefore,
the algebra $\mI_n$ is not a left coherent algebra, by Proposition
13.3, \cite{Stenstrom-RingQuot}. $\Box $

\begin{lemma}\label{a30Oct10}
$\ker_{\mI_2}(\cdot (H_1-H_2)) = \ker_{F_2}(\cdot (H_1-H_2)) =
\bigoplus_{i,j,k\in \N} Ke_{ij}(1) e_{kj}(2) \simeq
({}_{\mI_2}P_2)^{(\N )}$.
\end{lemma}

{\it Proof}. The algebra $B_2=\mI_2/ \ga_2$ is a domain
\cite{Bav-algintdif} where $\ga_2:= F(1) \t \mI_1(2) +\mI_1(1)\t
F(2)$  and $H_1-H_2\not\in \ga_2$. Therefore, $\CK := \ker_{\mI_2}
(\cdot (H_1-H_2))= \ker_{\ga_2} (\cdot (H_1-H_2))$. Let $F_2:=
F(1)\t F(2)$.  Notice that
$$ {}_{\mI_2}(\ga_2/F_2)_{\mI_2} \simeq F(1) \t B_1(2) \bigoplus
B_1(1)\t F(2)$$ is a direct sum of two $\mI_2$-bimodules. It
follows from the presentation $F(1) \t B_1(2) = \bigoplus_{i,j\in
\N, k\in \Z} e_{ij}(1) \t \der_2^k K[H_2]$ that $\ker_{F(1) \t
B_1(2)}(\cdot (H_1-H_2))=0$. Similarly, $\ker_{B_1(1)\t F(2)}
(\cdot (H_1-H_2)) =0$ (or use the $(1,2)$-symmetry). Therefore,
$$\CK = \ker_{F_2} (\cdot (H_1-H_2)) =\bigoplus_{i,j,k\in \N}
Ke_{ij}(1) e_{kj}(2)= \bigoplus_{j\in \N} (\bigoplus_{i,k\in \N}
Ke_{ij}(1) e_{kj}(2))\simeq \bigoplus_{j\in \N} ({}_{\mI_2}
P_2)\simeq {}_{\mI_2}(P_2)^{(\N )}.\;  \Box $$



\section{The algebra $\mI_n$ is maximal order}\label{AIMAXO}

The aim of this section is to prove Theorem \ref{14Feb11}.

Let $R$ be a ring. An element $x\in R$ is {\em right regular} if
$xr=0$ implies $r=0$ for $r\in R$. Similarly, a {\em left regular
element} is defined. A left and right regular element is called a
{\em regular element}. The sets of regular/left regular/right
regular elements of a ring $R$ are denoted respectively by
$\CC_R(0)$, $'\CC_R(0)$ and $\CC_R'(0)$.  For an arbitrary ring
$R$ there exists the {\em largest} (w.r.t. inclusion) {\em left
regular denominator set} $S_{l,0}=S_{l,0}(R)$ in the ring $R$
(regular means that $S_{l,0}(R)\subseteq \CC_R(0)$), and so
$Q_l(R):= S_{l,0}^{-1}R$ {\em is the largest left quotient ring
of} $R$ (Theorem 2.1, \cite{Bav-largquot}). Similarly, for an
arbitrary ring $R$ there exists the {\em largest right regular
denominator set} $S_{r,0}=S_{r,0}(R)$ in $R$, and so $Q_r(R):=
RS_{r,0}^{-1}$ is the {\em largest right quotient ring} of $R$.
The rings $Q_l(R)$ and $Q_r(R)$ were introduced and  studied in
\cite{Bav-largquot}. Let $\End_K(K[x])$ be the algebra of all
linear maps from the vector space $K[x]$ to itself and
$\Aut_K(K[x])$ be its group of units (i.e. the group of all
invertible linear maps from $K[x]$ to itself). The algebra $\mI_1$
is a subalgebra of $\End_K(K[x])$. Theorem 5.6.(1),
\cite{Bav-largquot}, states that $S_{r,0}(\mI_1)= \mI_1\cap
\Aut_K(K[x])$, it is the set of all elements of the algebra
$\mI_1$ that are invertible linear maps in $K[x]$. The set
$S_{r,0}(\mI_1)$ is huge comparing to the group of units $\mI_1^*$
of the algebra $\mI_1$ which is obviously a subset of
$S_{r,0}(\mI_1)$.

Let $R$ be a ring. A subring $S$ (not necessarily with 1)  of the
largest right quotient ring $Q_r(R)$ of the ring $R$ is called a
{\em right order} in $Q_r(R)$ if each element $q\in Q_r(R)$ has
the form $rs^{-1}$ for some elements $r,s\in S$. A subring $S$
(not necessarily with 1) of the largest left quotient ring
$Q_l(R)$ of the ring $R$ is called a {\em left order} in $Q_l(R)$
if each element $q\in Q_l(R)$ has the form $s^{-1}r$ for some
elements $r,s\in S$.

Let  $R_1$ and $R_2$ be right orders in $Q_r(\mI_n)$. We say that
the right orders $R_1$ and $R_2$ are {\em equivalent}, $R_1\sim
R_2$, if there are units $a_1, a_2, b_1, b_2\in Q_r(\mI_n)$ such
that $a_1R_1b_1\subseteq R_2$ and $a_2R_2b_2\subseteq R_1$.
Clearly, $\sim$ is an equivalent relation on the set of right
orders in $Q_r(\mI_n)$. A right order in $Q_r(\mI_n)$ is called a
{\em maximal right order} if it is maximal (w.r.t. $\subseteq $)
within its equivalence class.

\begin{lemma}\label{d12Feb11}
Let $Q_r(\mI_n)$ be the  right quotient ring of $\mI_n$ and $R, S$
be equivalent right orders in $Q_r(\mI_n)$ such that $R\subseteq
S$. Then there are equivalent right orders $T$ and $T'$ in
$Q_r(\mI_n)$ with $R\subseteq T\subseteq S$, $R\subseteq
T'\subseteq S$ and units $r_1, r_2$ of $Q_r(\mI_n)$ contained in
$R$ such that $r_1S\subseteq T$, $Tr_2\subseteq R$ and
$Sr_2\subseteq T'$, $r_1T'\subseteq R$. In particular,
$r_1Sr_2\subseteq R$.
\end{lemma}

{\it Proof}. By definition, $aSb\subseteq R$ for some units $a,b$
of $Q_r(\mI_n)$. Then $a= r_1s_1^{-1}$ and $b=r_2s_2^{-1}$ with
$r_i, s_i\in R$. Then $r_1Sr_2\subseteq r_1s_1^{-1}Sr_2\subseteq
Rs_2\subseteq R$. It is readily checked that $T=R+r_1S+Rr_1S$ and
$T' = R+Sr_2+Sr_2R$ are as claimed. $\Box $

$\noindent $

\begin{lemma}\label{a14Feb11}

\begin{enumerate}
\item $\CC_{\mI_n}(0)\cap \ga_n = \emptyset$, $'\CC_{\mI_n}(0)\cap
\ga_n = \emptyset$ and $\CC'_{\mI_n}(0)\cap \ga_n =
\emptyset$.\item $S_{l,0}(\mI_n) \cap \ga_n = \emptyset$ and
$S_{r,0}(\mI_n) \cap \ga_n = \emptyset$. \item For all elements $
a\in S_{l,0}(\mI_n) \cup S_{r,0}(\mI_n)$, $\mI_n a \mI_n = \mI_n$.
\end{enumerate}
\end{lemma}

{\it Proof}. 1. Trivial (since every element of the ideal $\ga_n$
is a left and right zero divisor in $\mI_n$).

2. Statement 2 follows from statement 1 and the inclusions $
S_{l,0}(\mI_n),S_{r,0}(\mI_n) \subseteq \CC_{\mI_n}(0)$.

3. If  $\mI_n a \mI_n \neq  \mI_n$ for some element $ a\in
S_{l,0}(\mI_n) \cup S_{r,0}(\mI_n)$ then $a\in \mI_n
a\mI_n\subseteq \ga_n$ (as $\ga_n$ is the only maximal ideal of
the algebra $\mI_n$). This contradicts to statement 2. $\Box $

$\noindent $

\begin{theorem}\label{14Feb11}
The algebra $\mI_n$ is a maximal left order in $Q_l(\mI_n)$ and a
maximal right order in $Q_r(\mI_n)$.
\end{theorem}

{\it Proof}. Suppose that $\mI_n\subseteq S$ and $S\sim \mI_n$ for
some right order $S$ in $Q_r(\mI_n)$. Then $aSb\subseteq \mI_n$
for some elements $a,b\in \mI_n\cap Q_r(\mI_n)^*$, by Lemma
\ref{d12Feb11}, where $Q_r(\mI_n)^*$ is the group of units of the
algebra $Q_r(\mI_n)$. By Theorem 2.8, \cite{Bav-largquot},
$\mI_n\cap Q_r(\mI_n)^* = S_{r,0}(\mI_n)$. Then, by Corollary
\ref{a14Feb11}.(3), $\mI_n\supseteq \mI_n a Sb\mI_n= (\mI_n
a\mI_n) S (\mI_n b \mI_n) = \mI_n S\mI_n = S$, i.e. $\mI_n = S$.
 Then the algebra $\mI_n$ is a maximal right order in
 $Q_r(\mI_n)$. Since the algebra $\mI_n$ admits an involution
 \cite{Bav-algintdif}, the algebra $\mI_n$ is also a maximal left
 order in $Q_l(\mI_n)$.  $\Box $


\small{

Department of Pure Mathematics

University of Sheffield

Hicks Building

Sheffield S3 7RH

UK

email: v.bavula@sheffield.ac.uk }


\begin{thebibliography}{99}
\bibitem{Bav-cdixintdif} V. V. Bavula, An analogue of the Conjecture of Dixmier is true for the algebra
 of polynomial integro-differential operators, (Arxiv:math.RA: 1011.3009).
\end{thebibliography}
\end{document}